\documentclass[11pt,reqno]{amsart}
\usepackage[a4paper,top=25mm,bottom=25mm,left=25mm,right=25mm]{geometry}
\usepackage[T1]{fontenc}
\usepackage{lmodern,mathtools,amssymb,microtype}
\usepackage{xcolor}
\usepackage[colorlinks=true,linkcolor=blue!45!black,citecolor=blue!45!black,urlcolor=blue!45!black]{hyperref}
\hypersetup{pdftitle={Riesz transforms on complete Riemannian manifolds for p > 2},pdfauthor={Peng Chen, Xuan Thinh Duong, Ji Li, and Lixin Yan}}
\newcommand{\dd}{\nabla}
\newcommand{\dif}{\mathrm d}
\newcommand{\one}{\mathbf 1}
\newcommand{\norm}[1]{\left\lVert#1\right\rVert}

\newtheorem{theorem}{Theorem}
\newtheorem{lemma}{Lemma}[section]
\numberwithin{equation}{section}
\allowdisplaybreaks[2]
\DeclareRobustCommand{\highprange}{\ensuremath{p>2}}
\title[Riesz transforms on complete Riemannian manifolds for \highprange]{Riesz transforms on complete Riemannian\\ manifolds for \highprange}
\author{Peng Chen}
\address{Peng Chen, Department of Mathematics, Sun Yat-sen University,
Guangzhou 510275, China}
\email{chenpeng3@mail.sysu.edu.cn}

\author{Xuan Thinh Duong}
\address{X. Duong, School of Mathematical and Physical Sciences,
Macquarie University, NSW 2109, Australia}
\email{xuan.duong@mq.edu.au}

\author{Ji Li}
\address{J. Li, School of Mathematical and Physical Sciences,
Macquarie University, NSW 2109, Australia}
\email{ji.li@mq.edu.au}

\author{Lixin Yan}
\address{Lixin Yan, Department of Mathematics, Sun Yat-sen University,
Guangzhou 510275, China}
\email{mcsylx@mail.sysu.edu.cn}

\date{}
\begin{document}
\begin{abstract}
Let $\Delta$ be the nonnegative Laplace--Beltrami operator on a complete
Riemannian manifold $M$ without boundary. Suppose that its heat semigroup satisfies
\[
 |\dd e^{-t\Delta}f|^2\le C_0 e^{-t\Delta}(|\dd f|^2),
 \qquad f\in C_c^\infty(M),\quad t>0.
\]
We prove that the Riesz transform $\dd\Delta^{-1/2}$ is bounded on $L^p(M)$
for every $2<p<\infty$, with norm depending only on $p$ and $C_0$.
This provides an affirmative answer to an open question raised by Coulhon and Duong \cite[Conjecture~1.1]{CD}.
\end{abstract}
\maketitle

\section{Introduction}
Let $M$ be a complete Riemannian manifold without boundary, with Riemannian
measure $\mu$ and gradient $\dd$. We write $|\cdot|$ for the length of a
tangent vector and $\|\cdot\|_p$ for the norm in $L^p(M,\mu)$,
$1\le p\le\infty$. Let $\Delta$ be the nonnegative  Laplace-Beltrami operator on $M$, and
$H_t=e^{-t\Delta}$ be the associated heat semigroup. The Riesz transform takes
functions to vector fields and has the representation
\begin{equation}\label{eq:Rp}
 \dd\Delta^{-1/2}f
 =\frac{1}{\sqrt{\pi}}\int_0^\infty \dd H_tf\,\frac{\dif t}{\sqrt t}.
\end{equation}
Strichartz \cite{Strichartz} asked for which $1<p<\infty$ the Riesz
transform satisfies
\begin{equation}\label{eq:R}
 \|\dd\Delta^{-1/2}f\|_p\le C_p\|f\|_p,\qquad f\in C_c^\infty(M).
\end{equation}
Integration by parts shows that
\[
\| |\dd f| \|_2=  \|\Delta^{1/2}f\|_2, \ \ \ \forall f\in C_c^\infty(M),
\]
therefore the Riesz transform  $\dd\Delta^{-1/2}$ is   bounded  from $L^2$ into itself.
For $p\ne2$, the estimate
\eqref{eq:R} has been studied on many classes of manifolds. See, for
example, \cite{A,AC,ACDH,B1985,B1987,B1989,CCH,CD1999,CD,CJKS,JL,L,S2004,Strichartz}
and the references therein.

In 2003, Coulhon and the second author \cite[Conjecture~1.1]{CD} conjectured that
\[
 \|\dd f\|_p\le C_p\|\sqrt{\Delta}f\|_p,\qquad f\in C_c^\infty(M),
\]
holds
for $2<p<\infty$ on every such manifold whose heat semigroup satisfies
\begin{equation}\label{eq:GE}
 |\dd H_tf|^2\le C_0H_t(|\dd f|^2),\qquad f\in C_c^\infty(M),\quad t>0.
\end{equation}
Here $C_0\ge1$ is independent of $f$ and $t$. This is condition (1.5)
of \cite{CD}.

We prove this assertion
for $2<p<\infty$ under \eqref{eq:GE}.

\begin{theorem}\label{thm:main}
Let $M$ be a complete Riemannian manifold without boundary satisfying
\eqref{eq:GE}. For every $2< p<\infty$, the Riesz transform
$\dd\Delta^{-1/2}$ extends to a bounded operator from $L^p(M)$ to
$L^p(TM)$ with norm at most $C_{p,C_0}$, depending only on $p$ and $C_0$.
In particular,
\[
 \|\dd f\|_p\le C_{p,C_0}\|\sqrt{\Delta}f\|_p,\qquad f\in C_c^\infty(M).
\]
\end{theorem}

\medskip
\noindent\textit{Remark.}
For $1<p<2$, Chen, Jiang, B. Li and H.-Q. Li \cite[Corollary~1.2]{CJLL} proved
the assertion  $\|\dd f\|_p\le C_p\|\sqrt{\Delta}f\|_p$ for $1<p<2$ via establishing the weak type $(1,1)$ estimate.

For $1<p<\infty$ and $q=p/(p-1)$, the estimate
$\|\dd f\|_p\le C\|\sqrt{\Delta}f\|_p$ for $f\in C_c^\infty(M)$ implies
by duality that
\[
 \|\sqrt{\Delta}f\|_q\le C\|\dd f\|_q,\qquad f\in C_c^\infty(M).
\]
See \cite[p.~172]{B1985}, \cite[p.~36]{B1989} or \cite[Proposition~2.1]{CD}.

\medskip
We compare our Theorem~\ref{thm:main} with earlier results. Write $B(x,r)$
for the open ball of radius $r$ centered at $x$, $V(x,r)=\mu(B(x,r))$
for its volume, and $h_t(x,y)$ for the heat kernel of $H_t$.
For manifolds with nonnegative Ricci curvature, Bakry \cite{B1987} used
Littlewood--Paley theory to prove \eqref{eq:R} for all $1<p<\infty$.
See also \cite{C1987,L1991}. A local change of metric can destroy the
nonnegativity of the Ricci curvature. Coulhon and Duong \cite{CD1999}
proved that the Riesz transform is of weak type $(1,1)$, and hence
bounded on $L^p$ for $1<p<2$, under the volume doubling condition
and the Gaussian heat kernel upper bound.

These two conditions are stable under uniformly equivalent Riemannian
metrics \cite[\S1.4]{ACDH}. Coulhon and Duong also gave an example
showing that these conditions do not suffice for the full range $p>2$
\cite{CD1999}.

Auscher, Coulhon, Duong and Hofmann \cite[Theorem~1.3]{ACDH} considered
complete noncompact manifolds satisfying the volume doubling condition and the Poincar\'e
inequality
\begin{equation}\label{eq:P}
 \int_B|f-f_B|^2\,\dif\mu\le Cr^2\int_B|\dd f|^2\,\dif\mu
\end{equation}
for every ball $B=B(x,r)$ and every smooth $f$, where
$f_B=\mu(B)^{-1}\int_B f\,\dif\mu$. For $2<p_0<\infty$, they proved
that the following assertions are equivalent.
\begin{itemize}
\item[(i)] For every $2<p<p_0$, there is a constant $C_p$ such that
$\|\dd H_tf\|_p\le C_pt^{-1/2}\|f\|_p$ for all $t>0$ and $f\in L^p(M)$.
\item[(ii)] The Riesz transform $\dd\Delta^{-1/2}$ is bounded on $L^p(M)$
for every $2<p<p_0$.
\end{itemize}
They also proved that, under the volume doubling condition
and the Gaussian heat kernel upper bound, the
pointwise gradient estimate
\begin{equation}\label{eq:G}
 |\dd_x h_t(x,y)|\le \frac{C}{\sqrt t\,V(x,\sqrt t)}
 \exp\!\left(-c\,\frac{d(x,y)^2}{t}\right),\qquad x,y\in M,\quad t>0,
\end{equation}
implies \eqref{eq:R} for every $1<p<\infty$ \cite[Theorem~1.4]{ACDH}.
For methods of proving heat kernel gradient estimates, see
\cite{LY,ACDH,CJKS}. Under the volume doubling condition, the Gaussian heat kernel upper bound and \eqref{eq:P},
the bound \eqref{eq:G} gives
\[
 |\dd H_tf|^2\le C H_{ct}(|\dd f|^2),\qquad f\in C_c^\infty(M),\quad t>0,
\]
for some $C,c>0$ \cite[Lemma~3.3]{ACDH}. Condition \eqref{eq:GE}
uses the same time $t$ on both sides.

Under \eqref{eq:GE}, Coulhon and Duong \cite[Theorem~4.1]{CD} proved
that for $2<p<\infty$,
\begin{equation}\label{eq:cdR}
 \|\dd f\|_p^2\le C_{p,C_0}\|f\|_p\|\Delta f\|_p,\qquad f\in C_c^\infty(M).
\end{equation}
On every complete Riemannian manifold, their Proposition~2.2 gives
\begin{equation}\label{eq:cd-sqrt}
 \|\sqrt{\Delta}f\|_p^2\le C_p\|f\|_p\|\Delta f\|_p,\qquad f\in C_c^\infty(M).
\end{equation}
These estimates give a common upper bound for the two squared norms.
The question was whether \eqref{eq:GE} alone implies
$\|\dd f\|_p\le C_{p,C_0}\|\sqrt{\Delta}f\|_p$.

Theorem~\ref{thm:main} answers this question by proving the Riesz
transform estimate \eqref{eq:R} for every $2<p<\infty$ under
\eqref{eq:GE} alone. We assume no volume doubling, curvature bound,
pointwise heat kernel upper bound, or Poincar\'e inequality~\eqref{eq:P}.

\section{A maximal estimate for Poisson averages of squares}\label{sec:averages}

\subsection{Scalar semigroup bounds}
The heat kernel is nonnegative and symmetric, with total mass at most
one \cite[Theorems 3.5 and 3.6]{Strichartz}. Define the Poisson semigroup $P_t=e^{-t\sqrt \Delta}$.  The Poisson subordination formula
\cite[\S1.3(ii), (1.1)]{CDsquare} reads
\begin{equation}\label{eq:subord}
 P_t f=\int_0^\infty H_s f\,\eta_t(s)\,\dif s,
\end{equation}
where
$$ \eta_t(s)=\frac{t}{2\sqrt\pi\,s^{3/2}}e^{-t^2/(4s)}.
$$
We see that $\int_0^\infty\eta_t(s)\,\dif s=1 .$

Thus both heat semigroups and Poisson semigroups have positive symmetric kernels of mass at most one.
Jensen's inequality and symmetry give contractivity on every $L^r$,
$1\le r\le\infty$. For $1<r<\infty$, the scalar maximal and square-function
estimates are
\begin{align}
 \norm{\sup_{s>0}|H_s f|}_r+\norm{\sup_{t>0}|P_t f|}_r
    &\le C_r\norm{f}_r,\label{eq:maximal}\\
 \norm{\left(\int_0^\infty |t\sqrt{\Delta} P_t f|^2\,\frac{\dif t}{t}\right)^{1/2}}_r
    &\le C_r\norm{f}_r.\label{eq:timeg}
\end{align}
For the proof, see Coulhon, Duong and Li \cite[p.~43, Lemma~2.1 and p.~39, \S1.3]{CDsquare}. See also Xu \cite[p.~3234, Theorem~8.1(i)]{Xu}.

\subsection{Spatial derivatives of Poisson averages}
We first provide the pointwise gradient estimates for the Poisson semigroup.
\begin{lemma}
\label{lem:poisson-gradient}
Let $M$ be a complete Riemannian manifold without boundary whose heat
semigroup satisfies \eqref{eq:GE} with constant $C_0$. For every $t>0$,
the following bounds for $P_t=e^{-t\sqrt{\Delta}}$ hold almost everywhere on $M$.
\begin{align}
 t^2|\dd P_t w|^2 &\le \frac{4\sqrt{2}}{e}\,C_0 P_{t/\sqrt2}(|w|^2)
       &&(w\in L^2(M)),\label{eq:poissmooth}\\*
 |\dd P_t u|^2 &\le C_0 P_t(|\dd u|^2)
       &&(u\in W^{1,2}(M)).\label{eq:dompoisson}
\end{align}
\end{lemma}

\begin{proof}
Completeness makes $C_c^\infty(M)$ a form core. Consequently
\eqref{eq:GE} extends to $W^{1,2}$.

For real $w,\phi\in C_c^\infty(M)$ with $\phi\ge0$, the heat equation and
the product rule for the energy form give, for $0<a<s$,
\begin{align*}
 \frac{\dif}{\dif a}\int_M (H_a\phi)(H_{s-a}w)^2\,\dif\mu
 &= -\int_M\!\langle\dd H_a\phi,\dd [(H_{s-a}w)^2]\rangle\,\dif\mu\\
 &\quad +2\int_M\!\langle\dd H_{s-a}w,\dd((H_a\phi)H_{s-a}w)\rangle\,\dif\mu\\
 &= 2\int_M (H_a\phi)|\dd H_{s-a}w|^2\,\dif\mu.
\end{align*}
Here and in the following sequel, $\langle \cdot,\cdot\rangle$ denotes the usual vector inner product.
The bounded form-domain factors admit the displayed product rule and the
endpoint values follow by strong $L^2$ continuity. Integrating in $a$ and
using symmetry of $H_a$ gives for all compact smooth $\phi$
\[
 \int_M (\phi)(H_s(w^2)-(H_s w)^2)\,\dif\mu
   =2\int_M (\phi) \left(\int_0^s H_a(|\dd H_{s-a}w|^2)\,\dif a\right)\,\dif\mu
\]
and thus almost everywhere we have
\[
 H_s(w^2)-(H_s w)^2
   =2\int_0^s H_a(|\dd H_{s-a}w|^2)\,\dif a
   \ge \frac{2s}{C_0}|\dd H_s w|^2 ,
\]
where in the last inequality we apply \eqref{eq:GE} to $H_{s-a}w\in W^{1,2}$.
Hence
$$
s|\dd H_s w|^2\le (C_0/2)H_s(|w|^2).
$$
Summing the estimates for real and imaginary parts gives the complex version.

Now subordination formula \eqref{eq:subord} and Jensen's
inequality for the probability measure $\eta_t(s)\dif s$  gives
\[
 t^2|\dd P_t w|^2
 \le \frac{C_0}{2}\int_0^\infty H_s(|w|^2)\frac{t^2}{s}\eta_t(s)\,\dif s
 \le \frac{4\sqrt2}{e}C_0\,P_{t/\sqrt2}(|w|^2).
\]
Indeed, with $z=t^2/s$ the ratio $(t^2/s)\eta_t(s)/\eta_{t/\sqrt2}(s)$
is $\sqrt2\,z e^{-z/8}\le 8\sqrt2/e$. This proves \eqref{eq:poissmooth}.

For \eqref{eq:dompoisson}, let $u\in W^{1,2}(M)$. subordination formula \eqref{eq:subord}, Jensen's inequality
and \eqref{eq:GE} give
\[
 |\dd P_t u|^2\le\int_0^\infty |\dd H_su|^2\eta_t(s)\,\dif s
 \le C_0\int_0^\infty H_s(|\dd u|^2)\eta_t(s)\,\dif s
 =C_0 P_t(|\dd u|^2).
\]
The proof is complete.
\end{proof}

For every finite family $w_j\in L^r\cap L^2$, $t_j>0$, we also have the following vector valued square function estimates
\begin{equation}\label{eq:vector}
 \norm{\bigg(\sum_j|t_j\dd P_{t_j}w_j|^2\bigg)^{1/2}}_r
 \le C_r\norm{\bigg(\sum_j |w_j|^2\bigg)^{1/2}}_r,
 \qquad 1<r\le2.
\end{equation}
Indeed, $\int_0^\infty t^2\eta_t(s)^2\,\dif s=\pi^{-1}$, so Cauchy--Schwarz
in \eqref{eq:subord} gives
\[
 \sum_j |t_j\dd P_{t_j}w_j|^2
 \le \pi^{-1}\sum_j\int_0^\infty|\dd H_sw_j|^2\,\dif s.
\]
The spatial heat square-function theorem
\cite[p.~41, Theorem~1.2]{CDsquare} yields
$$\norm{\bigg(\int_0^\infty|\dd H_s w|^2\,\dif s\bigg)^{1/2}}_r
 \le C_r\norm{w}_r$$ for $1<r\le2$.
Let $\varepsilon_j$ be independent random signs and $\mathbb E$ their
expectation. Apply the Hilbert-space Khintchine inequality in
$L^2((0,\infty),\dif s;T_xM)$ at each $x$, followed by the scalar
square-function bound for $\sum_j\varepsilon_j w_j$. This gives
\begin{align*}
 \int_M\left(\sum_j\int_0^\infty |\dd H_s w_j|^2\,\dif s\right)^{r/2}\dif\mu
 &\le C_r\mathbb E\int_M
  \left(\int_0^\infty\left|\dd H_s\Big(\sum_j\varepsilon_j w_j\Big)\right|^2
                    \,\dif s\right)^{r/2}\dif\mu\\
 &\le C_r\int_M \mathbb E\left|\sum_j\varepsilon_j w_j\right|^r\dif\mu
 \le C_r\int_M\left(\sum_j|w_j|^2\right)^{r/2}\dif\mu .
\end{align*}
The last inequality uses $r\le2$ and orthogonality of the signs.

\subsection{The weighted Rellich identity}

Take real $f\in C_c^\infty(M)$  and $w\in L^2(M)\cap L^\infty(M)$ with $w\ge0$, and fix $t>0$.
The following differentiation formulas hold for every $s>0$.
\begin{align}
 \partial_s\int_s^{s+t}P_u f\,\dif u
   &=-(P_s-P_{s+t})f,\notag\\
 \partial_s^2\int_s^{s+t}P_u f\,\dif u
   &=\sqrt{\Delta}(P_s-P_{s+t})f
    =\Delta\int_s^{s+t}P_u f\,\dif u,\notag\\
 \dd\int_s^{s+t}P_u f\,\dif u
   &=\dd\Delta^{-1/2}(P_s-P_{s+t})f.\label{eq:rel-diff}
\end{align}
By \eqref{eq:poissmooth},
$P_tw$ and $\dd P_tw$ are bounded. We may therefore use
$(P_tw)(P_s-P_{s+t})f$ as a form test function.
Differentiation followed by one integration by parts gives
\begin{align}
 &\frac12\frac{\dif}{\dif s}\int_M(P_tw)
 \left(\left|\dd\int_s^{s+t}P_u f\,\dif u\right|^2
              -|(P_s-P_{s+t})f|^2\right)\dif\mu\notag\\
 &\quad=-\int_M(P_tw)\left\langle
      \dd\int_s^{s+t}P_u f\,\dif u,\dd(P_s-P_{s+t})f\right\rangle\dif\mu\notag\\
 &\qquad +\int_M(P_tw)((P_s-P_{s+t})f)\,
        \Delta\!\int_s^{s+t}P_u f\,\dif u\,\dif\mu\notag\\
 &\quad=\int_M((P_s-P_{s+t})f)
          \left\langle\dd P_tw,\dd\int_s^{s+t}P_u f\,\dif u\right\rangle\dif\mu .
 \label{eq:rel-derivative}
\end{align}
At $s\downarrow0$, $(P_s-P_{s+t})f\to(I-P_t)f$ in $L^2$. At
$s\to\infty$ the limit is zero, since $(I-P_t)f\perp\ker(\sqrt{\Delta})$. By
\eqref{eq:rel-diff} the limits of the spatial derivatives follow by applying
the $L^2$ contraction $\dd\Delta^{-1/2}$. Integrating \eqref{eq:rel-derivative} gives
\begin{align}
 &\int_M(P_tw)|\dd\Delta^{-1/2}(I-P_t)f|^2\,\dif\mu\notag\\
 &=\int_M(P_tw)|(I-P_t)f|^2\,\dif\mu\notag\\
 &\quad-2\int_0^\infty\!\int_M ((P_s-P_{s+t})f)
       \left\langle\dd P_tw,\int_s^{s+t}\dd P_u f\,\dif u\right\rangle
       \dif\mu\,\dif s\notag\\
 &\eqqcolon \mathcal S_t(f,w)+\mathcal E_t(f,w).
 \label{eq:rellich}
\end{align}
Here $\mathcal S_t(f,w)$ is the scalar term and $\mathcal E_t(f,w)$ is
the remainder, including the factor $-2$. The estimates below give
absolute convergence of the double integral at the
dyadic times used in the proof. For each fixed $t>0$ the same estimates
apply to the single interval indexed by $t$.

\subsection{Summation over dyadic times}

In this subsection our aim is to prove, for $4\le p<\infty$,
\begin{equation}
 \left\|\sup_{i\in\mathbb Z}P_{2^i}
       \bigl(|\dd\Delta^{-1/2}(I-P_{2^i})f|^2\bigr)\right\|_{p/2}
 \le C_{p,C_0}\|f\|_p^2,\qquad f\in L^2\cap L^p .
 \label{eq:maxsq}
\end{equation}
Set $r=(p/2)'=p/(p-2)$ and $t_i=2^i$. For every finite nonempty
$J\subset\mathbb Z$ and every choice of nonnegative $w_i\in L^2(M)\cap L^\infty(M)$,
$i\in J$, we first prove the weighted estimate ,
\begin{align}
 \sum_{i\in J}\int_M (P_{t_i}w_i)\bigl(|R(I-P_{t_i})f|^2\bigr)\,\dif\mu
 \le C_{p,C_0}\|f\|_p^2\left\|\sum_{i\in J}w_i\right\|_r .
 \label{eq:rel-weighted1}
\end{align}
The constant is independent of $J$ and the weights.

From \eqref{eq:rellich}, we estimate $\sum_i\mathcal S_{t_i}(f,w_i)$ and
$\sum_i|\mathcal E_{t_i}(f,w_i)|$ separately.
For the scalar terms, symmetry, positivity, and the Poisson maximal inequality~\eqref{eq:maximal} give
\begin{align}
 \sum_i\mathcal S_{t_i}(f,w_i)
 &\le \int_M\Big(\sum_i w_i\Big)\sup_{a>0}P_a\!\left(\left[|f|+\sup_{u>0}|P_u f|\right]^2\right)\dif\mu\notag\\
 &\le C_p\left\|\sum_i w_i\right\|_r\|f\|_p^2.
 \label{eq:rel-scalar}
\end{align}

To estimate the remainders $\mathcal E_{t_i}(f,w_i)$, put $s=vt_i$.
Write $I_{i,v}=[vt_i,(v+1)t_i]$. Then
\[
 (P_{vt_i}-P_{(v+1)t_i})f=\int_{I_{i,v}} \sqrt{\Delta} P_u f\,\dif u,
 \qquad \int_{I_{i,v}}\frac{\dif u}{u}=\log(1+1/v).
\]
For the spatial factor, \eqref{eq:poissmooth} gives pointwise
\[
 \max_i\left|\int_{I_{i,v}}\dd P_u f\,\dif u\right|
 \le C\sqrt{C_0}\,\log(1+1/v)\left[\sup_{a>0}P_a(|f|^2)\right]^{1/2}.
\]
Consequently,
\begin{equation}
 \left\|\max_i\left|\int_{I_{i,v}}\dd P_u f\,\dif u\right|\right\|_p
 \le C_{p,C_0}\log(1+1/v)\|f\|_p .
 \label{eq:rel-spatial}
\end{equation}
For the scalar factor, Cauchy--Schwarz with measure $\dif u/u$ gives
\[
 \left|\int_{I_{i,v}} \sqrt{\Delta} P_u f\,\dif u\right|^2
 \le\log(1+1/v)\int_{I_{i,v}} |u \sqrt{\Delta} P_u f|^2\,\frac{\dif u}{u}.
\]
For fixed $u>0$ the condition $u\in I_{i,v}$ is equivalent to
$\log_2(u/(v+1))\le i\le\log_2(u/v)$.
The number of such integers is at most $1+\log_2(1+1/v)$.
Sum in $i$ and apply \eqref{eq:timeg} to get
\begin{equation}
 \left\|\left(\sum_i\left|\int_{I_{i,v}} \sqrt{\Delta} P_u f\,\dif u\right|^2
       \right)^{1/2}\right\|_p
 \le C_p\sqrt{\log(1+1/v)\Big(1+\log(1+1/v)\Big)}\,\|f\|_p .
 \label{eq:rel-time}
\end{equation}
For the remaining factor, \eqref{eq:vector} and $w_i\ge0$ give
\begin{equation}
 \left\|\left(\sum_i|t_i\dd P_{t_i} w_i|^2\right)^{1/2}\right\|_r
 \le C_r\left\|\left(\sum_i |w_i|^2\right)^{1/2}\right\|_r
 \le C_r\left\|\sum_i w_i\right\|_r .
  \label{eq:real-weight1}
\end{equation}
By \eqref{eq:rel-spatial}, \eqref{eq:rel-time} and \eqref{eq:real-weight1}, using Cauchy--Schwarz in $i$ and H\"older in $x$ with $1/p+1/p+1/r=1$,
\begin{align*}
 &\sum_i|\mathcal E_{t_i}(f,w_i)|\\
 &\quad\le 2\sum_i\int_0^\infty\!\int_M |(P_s-P_{s+t_i})f|\,|\dd P_{t_i}w_i|
       \left|\int_s^{s+t_i}\dd P_u f\,\dif u\right|\dif\mu\,\dif s\\
 &\quad\le 2\int_0^\infty
 \left\|\max_i\left|\int_{I_{i,v}}\dd P_u f\,\dif u\right|\right\|_p
 \left\|\left(\sum_i\left|\int_{I_{i,v}} \sqrt{\Delta} P_u f\,\dif u\right|^2
       \right)^{1/2}\right\|_p
 \left\|\left(\sum_i|t_i\dd P_{t_i}w_i|^2\right)^{1/2}\right\|_r\dif v\\
 &\quad\le C_{p,C_0}\|f\|_p^2\left\|\sum_i w_i\right\|_r
           \int_0^\infty\big(\log(1+1/v)\big)^{3/2}\big(1+\log(1+1/v)\big)^{1/2}\,\dif v .
\end{align*}
For $0<v\le1$ the integrand is at most $C(1+|\log v|)^2$, and for
$v\ge1$ it is at most $Cv^{-3/2}$. Using this bound for the remainders
and \eqref{eq:rel-scalar} for the scalar terms, we sum \eqref{eq:rellich}
in $i$ to achieve \eqref{eq:rel-weighted1}. By symmetry of $P_{t_i}$, \eqref{eq:rel-weighted1} implies
\begin{equation}
 \sum_i\int_M w_i P_{t_i}\bigl(|\dd\Delta^{-1/2}(I-P_{t_i})f|^2\bigr)\,\dif\mu
 \le C_{p,C_0}\|f\|_p^2\left\|\sum_i w_i\right\|_r .
 \label{eq:rel-weighted}
\end{equation}

For duality, take $0\le w\in C_c^\infty(M)$ and order the finite
nonempty set $J\subset\mathbb Z$ by the usual order on the integers.
For $i\in J$, define
\begin{align*}
 E_i={}&\bigg(\bigcap_{\substack{j\in J\\ j<i}}
 \left\{P_{t_i}(|\dd\Delta^{-1/2}(I-P_{t_i})f|^2)>P_{t_j}(|\dd\Delta^{-1/2}(I-P_{t_j})f|^2)\right\}\bigg)\\
 &{}\bigcap\bigg(\bigcap_{\substack{j\in J\\ j>i}}
 \left\{P_{t_i}(|\dd\Delta^{-1/2}(I-P_{t_i})f|^2)\ge P_{t_j}(|\dd\Delta^{-1/2}(I-P_{t_j})f|^2)\right\}\bigg).
\end{align*}
  Thus $M=\bigsqcup_{i\in J}E_i$.  Estimate
\eqref{eq:rel-weighted} applies to $w_i=w\one_{E_i}$.
 Thus
\begin{align*}
 \int_M w\max_{i\in J}P_{t_i}(|\dd\Delta^{-1/2}(I-P_{t_i})f|^2)\,\dif\mu
 &= \sum_i \int_{E_i} w \max_{i\in J}P_{t_i}(|\dd\Delta^{-1/2}(I-P_{t_i})f|^2)\,\dif\mu\\
 &=\sum_i \int_{E_i} w_i P_{t_i}(|\dd\Delta^{-1/2}(I-P_{t_i})f|^2)\,\dif\mu\\
 &\le C_{p,C_0}\|f\|_p^2\|w\|_r .
\end{align*}
Duality and monotone convergence on $J=[-N,N]\cap\mathbb Z$ give
\eqref{eq:maxsq}. The proof also applies directly to real
$f\in L^2\cap L^p$. The spectral identities \eqref{eq:rel-diff}
and their endpoint limits use only $L^2$, and the bounds in the dyadic
sum use only \eqref{eq:maximal}, \eqref{eq:timeg}, and \eqref{eq:poissmooth}
on their stated domains, together with the weight estimate \eqref{eq:vector}. For complex $f$, add the bounds for its
real and imaginary parts.

\section{Proof of Theorem~\ref{thm:main}}\label{sec:riesz}

We start with a scalar function $F\ge0$. On a finite range of dyadic
times, the proof of Lemma~\ref{lem:stopping} groups the points of the set
$E$ below according to the last $t$ for which $P_tF>F/(2K)$. At every
later dyadic time $s$ in this range, these points satisfy $F>KP_sF$.
For $F=|\dd u|^2$ and $K=4C_0$, \eqref{eq:dompoisson} gives
$|\dd u|\le2|\dd u-\dd P_su|$ on $\{F>KP_sF\}$.
This allows us to use \eqref{eq:maxsq} with $f=\sqrt{\Delta}\,u$.

\subsection{A scalar inequality for Poisson averages}\label{sec:stopping}

Assume that $M$ is connected. For $2\le q<\infty$ and $F\in L^q(M)$, set
\[
 \Pi F=
 \begin{cases}
  \displaystyle\mu(M)^{-1}\int_M F\,\dif\mu,& \mu(M)<\infty,\\
  0,&\mu(M)=\infty.
 \end{cases}
\]
On $L^2$ this is the spectral projection onto $\ker(\sqrt{\Delta})$. Indeed,
$\|\sqrt{\Delta} h\|_2=\|\dd h\|_2$ identifies the kernel with the square-integrable
constants. In finite volume $1\in W^{1,2}(M)$ and $\dd 1=0$.

Auscher, Coulhon, Duong and Hofmann \cite[Lemma~2.2]{ACDH} prove a
good-$\lambda$ inequality using averages over balls on a doubling space.
In the following lemma the average at time $t$ is
$P_t(F\one_{\{F>KP_tF\}})$. The proof uses only Poisson averages and
requires no volume doubling assumption.

\begin{lemma}\label{lem:stopping}
Let $K>1$, $2\le q<\infty$ and $0\le F\in L^q(M)$. Then
\begin{equation}\label{eq:stop}
 \|F\|_q\le C_{q,K}
 \left(\left\|\sup_{j\in\mathbb Z}
   P_{2^j}\!\left(F\one_{\{F>KP_{2^j}F\}}\right)\right\|_q
       +\|\Pi F\|_q\right).
\end{equation}
\end{lemma}

\begin{proof}
Set $L=2K$, $r=q/2$. Fix integers $m<n$ and define
\[
 S_i=\{F>KP_{2^i}F\},\qquad
 E=\{F>0,\ P_{2^m}F>F/L,\ P_{2^n}F\le F/L\}.
\]
For each $i=m,\ldots,n-1$, define
\[
 E_i=E\cap\{P_{2^i}F>F/L\}
       \cap\bigcap_{j=i+1}^{n}\{P_{2^j}F\le F/L\}.
\]
Thus $x\in E_i$ precisely when $x\in E$ and $i$ is the last index
in $\{m,\ldots,n\}$ at which $P_{2^i}F(x)>F(x)/L$. For each
$x\in E$ the inequality holds at $m$ and fails at $n$, so this last
index exists and belongs to $\{m,\ldots,n-1\}$. Its uniqueness gives
$E=\bigsqcup_{i=m}^{n-1}E_i$. For $x\in E_i$ and $i<j\le n$,
\[
 KP_{2^j}F(x)\le (K/L)F(x)=F(x)/2<F(x),
 \qquad\text{hence } E_i\subset S_j .
\]
On $E_i$ the inequality $P_{2^i}F>F/L$ and the identity $q=2r$ give
\[
 L^{-r}F^q=F^r(F/L)^r\le F^r(P_{2^i}F)^r .
\]
For the second pointwise bound, let $p_t^{\rm P}(x,y)$ denote the kernel
of $P_t$. For almost every $x$, the measure
$p_t^{\rm P}(x,y)\,\dif\mu(y)$ has mass at most one. Add its missing mass
at a point where $F$ is assigned the value zero. Jensen's inequality
for the resulting probability measure and the convex function
$a\mapsto a^r$ on $[0,\infty)$ gives
\[
 (P_t F(x))^r
 =\left(\int_M F(y)p_t^{\rm P}(x,y)\,\dif\mu(y)\right)^r
 \le \int_M F(y)^r p_t^{\rm P}(x,y)\,\dif\mu(y)
 =P_t(F^r)(x).
\]
For $r=1$ this is an equality. Since
$\|F^r\one_{E_i}\|_2^2\le\|F^r\|_2^2=\int_M F^q\,\dif\mu<\infty$,
both functions in the following $L^2$ pairing are admissible.
Self-adjointness and $L^2$ contractivity give
\begin{align*}
 \int_{E_i}F^rP_{2^i}(F^r)\,\dif\mu
 &=\langle F^r\one_{E_i},P_{2^i}(F^r)\rangle_{L^2}
  =\langle P_{2^i}(F^r\one_{E_i}),F^r\rangle_{L^2},\\
 \big|\langle P_{2^i}(F^r\one_{E_i}),F^r\rangle_{L^2}\big|
 &\le\|F^r\one_{E_i}\|_2\,\|F^r\|_2 .
\end{align*}
Integrating over each $E_i$ and summing over the finite disjoint
partition $E=\bigsqcup_{i=m}^{n-1}E_i$ now yields
\begin{align}
 L^{-r}\int_E F^q\,\dif\mu
 &=\sum_{i=m}^{n-1}\int_{E_i}F^r(F/L)^r\,\dif\mu \notag\\
 &\le\sum_{i=m}^{n-1}\int_{E_i}F^r(P_{2^i}F)^r\,\dif\mu
 \le\sum_{i=m}^{n-1}\int_{E_i}F^rP_{2^i}(F^r)\,\dif\mu \notag\\
 &=\int_M F^r\sum_{i=m}^{n-1}P_{2^i}(F^r\one_{E_i})\,\dif\mu .
 \label{eq:stoppair}
\end{align}

Expanding the squared $L^2$ norm in \eqref{eq:stoppair} and using
$P_sP_t=P_{s+t}$ gives
\begin{align*}
 \left\|\sum_{i=m}^{n-1}P_{2^i}(F^r\one_{E_i})\right\|_2^2
 ={}&\sum_{i=m}^{n-1}\int_{E_i} F^r P_{2^{i+1}}(F^r\one_{E_i})\,\dif\mu
 +2\sum_{m\le j<i<n}\int_{E_i} F^r
                 P_{2^i+2^j}(F^r\one_{E_j})\,\dif\mu .
\end{align*}
In the first sum $E_i\subset S_{i+1}$. For the second sum, fix
$m<i<n$. We compare all its summands at time $2^i$ using
$P_bh\le(b/a)P_ah$ for $h\ge0$ and $0<a\le b$; see
\cite[Lemma~1.2]{Mei}. For convenience, we give the calculation here.
For $0<a\le b$ and $s>0$, the density in \eqref{eq:subord} satisfies
\[
 \frac{\eta_b(s)}{\eta_a(s)}
 =\frac ba \exp\!\left(-\frac{b^2-a^2}{4s}\right)\le\frac ba .
\]
Fix $m\le j<i$. Since $F^r\one_{E_j}\in L^2$ is nonnegative,
$H_s(F^r\one_{E_j})\ge0$. Subordination therefore gives, almost everywhere,
\begin{align*}
 P_b(F^r\one_{E_j})
 &=\int_0^\infty\eta_b(s)H_s(F^r\one_{E_j})\,\dif s
 \le\frac ba\int_0^\infty\eta_a(s)H_s(F^r\one_{E_j})\,\dif s
 =\frac ba P_a(F^r\one_{E_j}).
\end{align*}
Take $a=2^i$ and $b=2^i+2^j$. Since $b/a=1+2^{j-i}\le 3/2\le 2$, this gives
\[
 P_{2^i+2^j}(F^r\one_{E_j})\le 2P_{2^i}(F^r\one_{E_j}).
\]
Next, if $y\in E_j$ with $j<i$, the definition of $E_j$ gives
$P_{2^i}F(y)\le F(y)/L$. Since $L=2K$ and $F(y)>0$ on $E_j$,
\[
 KP_{2^i}F(y)\le F(y)/2<F(y),\qquad y\in S_i .
\]
Thus $\bigcup_{j=m}^{i-1}E_j\subset S_i$. Disjointness gives
$\sum_{j=m}^{i-1}\one_{E_j}=\one_{\cup_{j=m}^{i-1}E_j}\le\one_{S_i}$.
Summing the time comparisons, then using linearity and positivity, yields
\begin{align*}
 \sum_{j=m}^{i-1}P_{2^i+2^j}(F^r\one_{E_j})
 &\le 2\sum_{j=m}^{i-1}P_{2^i}(F^r\one_{E_j})
 =2P_{2^i}\!\left(F^r\sum_{j=m}^{i-1}\one_{E_j}\right)
 \le 2P_{2^i}(F^r\one_{S_i}).
\end{align*}
Next we estimate the term
$$
\left\|\sum_{i=m}^{n-1}P_{2^i}(F^r\one_{E_i})\right\|_2^2=\sum_{i=m}^{n-1}\sum_{j=m}^{n-1}\int_M P_{2^i}(F^r\one_{E_i}) P_{2^j}(F^r\one_{E_j})\,\dif\mu
$$
through diagonal part and off-diagonal part. Fix $m\le i<n$. The diagonal term in the squared norm has Poisson time
$2^i+2^i=2^{i+1}$. Since $E_i\subset S_{i+1}$, positivity gives
\begin{align*}
 \int_{E_i}F^r P_{2^{i+1}}(F^r\one_{E_i})\,\dif\mu
 &\le \int_{E_i}F^r P_{2^{i+1}}(F^r\one_{S_{i+1}})\,\dif\mu\\
 &\le \int_{E_i}F^r\sup_{k\in\mathbb Z}P_{2^k}(F^r\one_{S_k})\,\dif\mu .
\end{align*}
The last inequality uses the term $k=i+1$ in the supremum. At the
endpoint $i=n-1$, the inclusion $E_{n-1}\subset S_n$ follows from
$K P_{2^n}F\le F/2<F$ on $E_{n-1}\subset E$.
In the off-diagonal part, the pairs $(i,j)$ and $(j,i)$ give the same
inner product. Thus the terms with $j<i$ carry a factor of $2$.
Multiplying the preceding time comparison for the sum over $j<i$ by
$F^r\ge0$ and integrating over $E_i$ yields
\begin{align*}
 2\sum_{j=m}^{i-1}\int_{E_i}F^r P_{2^i+2^j}(F^r\one_{E_j})\,\dif\mu
 &=2\int_{E_i}F^r\sum_{j=m}^{i-1}P_{2^i+2^j}(F^r\one_{E_j})\,\dif\mu\\
 &\le 2\cdot2\int_{E_i}F^r P_{2^i}(F^r\one_{S_i})\,\dif\mu\\
 &\le 4\int_{E_i}F^r\sup_{k\in\mathbb Z}P_{2^k}(F^r\one_{S_k})\,\dif\mu .
\end{align*}
The second factor of $2$ comes from comparing $P_{2^i+2^j}$ with
$P_{2^i}$, and the last inequality uses $k=i$. For $i=m$ the sum over
$j<i$ is empty. Adding the two bounds for each $i$ and using
$\sum_{i=m}^{n-1}\one_{E_i}=\one_E$ gives
\begin{align}
 \left\|\sum_{i=m}^{n-1}P_{2^i}(F^r\one_{E_i})\right\|_2^2
 &\le (1+4)\sum_{i=m}^{n-1}\int_{E_i}F^r
       \sup_{k\in\mathbb Z}P_{2^k}(F^r\one_{S_k})\,\dif\mu \notag\\
 &=5\int_E F^r\sup_{k\in\mathbb Z}P_{2^k}(F^r\one_{S_k})\,\dif\mu .
 \label{eq:stopcross}
\end{align}

For $q>2$ put $\ell=(q+r)/2$ and $\theta=(\ell-r)/(\ell-1)$.
Then $1<r<\ell<q$, $0<\theta<1$, and $r=\theta+\ell(1-\theta)$.
H\"older's inequality for the Poisson kernel gives
\[
 P_{2^j}(F^r\one_{S_j})\le
    \big[P_{2^j}(F\one_{S_j})\big]^\theta
    \big[P_{2^j}(F^\ell)\big]^{1-\theta}.
\]
Taking suprema and applying H\"older on $M$, with
$\theta/q+(1-\theta)/(q/\ell)=1/2$, yields
\begin{align}
 \left\|\sup_j P_{2^j}(F^r\one_{S_j})\right\|_2
 &\le
 \left\|\sup_j P_{2^j}(F\one_{S_j})\right\|_q^\theta
 \left\|\sup_{t>0}P_t(F^\ell)\right\|_{q/\ell}^{1-\theta}
 \notag\\
 &\le C_q\left\|\sup_j P_{2^j}(F\one_{S_j})\right\|_q^\theta
        \|F\|_q^{r-\theta}.
 \label{eq:stopmax}
\end{align}
The last line uses \eqref{eq:maximal} at $q/\ell>1$.
For $q=2$, take $r=\theta=1$ and the last bound is an equality with
constant one.

Put $G=\sum_{i=m}^{n-1}P_{2^i}(F^r\one_{E_i})$ and
$U=\sup_j P_{2^j}(F^r\one_{S_j})$.
Since $2r=q$, we have $\|F^r\|_2=\|F\|_q^r$.
Each term of the finite sum $G$ belongs to $L^2$ by contractivity.
Also, $\sup_j P_{2^j}(F\one_{S_j})\le\sup_{t>0}P_tF\in L^q$ by
\eqref{eq:maximal}, so \eqref{eq:stopmax} gives $U\in L^2$.
By \eqref{eq:stoppair} and Cauchy--Schwarz,
\[
 \int_E F^q\,\dif\mu
 \le L^r\int_M F^rG\,\dif\mu
 \le L^r\|F\|_q^r\|G\|_2 .
\]
By \eqref{eq:stopcross} and Cauchy--Schwarz on $E$,
\[
 \|G\|_2^2\le 5\int_E F^r U\,\dif\mu
 \le 5\|F^r\one_E\|_2\|U\|_2\le 5\|F\|_q^r\|U\|_2 .
\]
Taking the square root and substituting in the first bound gives
$$\int_E F^q\,\dif\mu\le\sqrt{5}\,L^r\|F\|_q^{3r/2}\|U\|_2^{1/2}.$$
By \eqref{eq:stopmax},
\[
 \|U\|_2^{1/2}\le C_q
 \left\|\sup_jP_{2^j}(F\one_{S_j})\right\|_q^{\theta/2}
 \|F\|_q^{(r-\theta)/2}.
\]
Since $3r/2+(r-\theta)/2=2r-\theta/2=q-\theta/2$, we obtain
\begin{equation}\label{eq:stopE}
 \int_E F^q\,\dif\mu
 \le C_{q,K}\,\|F\|_q^{q-\theta/2}
       \left\|\sup_jP_{2^j}(F\one_{S_j})\right\|_q^{\theta/2}.
\end{equation}

For $b\in C_c^\infty(M)$, the spectral theorem gives $P_tb\to b$
as $t\downarrow0$ and $P_tb\to\Pi b$ as $t\to\infty$ in $L^2$.
Both differences have $L^\infty$ norm at most $2\|b\|_\infty$.
The inequality $\|h\|_q^q\le\|h\|_\infty^{q-2}\|h\|_2^2$ and
$L^q$ contractivity give both limits for test functions in $L^q$, and density extends them to all $F\in L^q$.
Since $F\ge0$, the definition of $E$ gives
\[
 E^c=\{F=0\}\cup\{P_{2^m}F\le F/L\}\cup\{P_{2^n}F>F/L\}.
\]
If $P_{2^m}F(x)\le F(x)/L$, then $L=2K>1$ gives
\[
 |F(x)-P_{2^m}F(x)|=F(x)-P_{2^m}F(x)
 \ge \frac{L-1}{L}F(x).
\]
If $P_{2^n}F(x)>F(x)/L$, then $F(x)<L P_{2^n}F(x)$.
Thus, at each point of $E^c$ where $F(x)>0$, at least one of
the two terms on the right below bounds $F(x)$. Both terms are
nonnegative by positivity of $P_{2^n}$. At points of $E$ or
of $\{F=0\}$ the left side vanishes. Hence, almost everywhere on $M$,
\[
 F\one_{E^c}\le\frac{L}{L-1}|F-P_{2^m}F|+L P_{2^n}F.
\]
Take $m=-N$, $n=N$, so that $E$ depends on $N$. The triangle
inequality and the two $L^q$ limits above give
\[
 \|F\one_{E^c}\|_q
 \le\frac{L}{L-1}\|F-P_{2^{-N}}F\|_q+L\|P_{2^N}F\|_q
 \xrightarrow[N\to\infty]{} L\|\Pi F\|_q .
\]
Thus $\limsup_{N\to\infty}\|F\one_{E^c}\|_q\le L\|\Pi F\|_q$.
Combining this with \eqref{eq:stopE} in the identity
$\|F\|_q^q=\int_E F^q\,\dif\mu+\|F\one_{E^c}\|_q^q$ gives
\[
 \|F\|_q^q\le C_{q,K}\|F\|_q^{q-\theta/2}
    \left\|\sup_j P_{2^j}(F\one_{S_j})\right\|_q^{\theta/2}
    +L^q\|\Pi F\|_q^q .
\]
Young's inequality absorbs $\tfrac12\|F\|_q^q$ and gives
\eqref{eq:stop}.
\end{proof}

\subsection{Recovering the spatial gradient}
Fix $4\le p<\infty$. We first prove
\begin{equation}\label{eq:form-to-R}
 \|\dd u\|_p\le C_{p,C_0}\|\sqrt{\Delta} u\|_p
\end{equation}
for $u\in W^{1,2}(M)$ with $\dd u\in L^p(TM)$ and $\sqrt{\Delta}\,u\in L^p(M)$.

For $t>0$, let $S_t=\{|\dd u|^2>4C_0P_t(|\dd u|^2)\}$.
Since $\dd u\in L^p(TM)$, we have $|\dd u|^2\in L^{p/2}(M)$.
By using \eqref{eq:stop} with $F=|\dd u|^2$, $q=p/2$ and $K=4C_0$, we see that
\begin{align}\label{dd u}
 \|\dd u\|_p^2\leq C_{p,C_0}\bigg(\left\|\sup_i P_{2^i}\bigl(|\dd u|^2\one_{S_{2^i}}\bigr)\right\|_{p/2}
+
 \|\Pi(|\dd u|^2)\|_{p/2}\bigg)=:Term_1+Term_2.
\end{align}

We first consider $Term_1$. 
Since $\Delta=\dd^*\dd$, the polar decomposition
of $\dd$ gives
\begin{equation}\label{eq:polar}
 (\dd\Delta^{-1/2})\sqrt{\Delta}\,u=\dd u
 \quad(u\in W^{1,2}(M)),\qquad
 \|\dd\Delta^{-1/2}f\|_2=\|(I-\Pi)f\|_2
 \quad(f\in L^2(M)).
\end{equation}
Then we write
\begin{align*}
 \dd u
 &=\dd\Delta^{-1/2}\sqrt{\Delta}\,u
 =\dd\Delta^{-1/2}P_t\sqrt{\Delta}\,u
   +\dd\Delta^{-1/2}(I-P_t)\sqrt{\Delta}\,u\\
 &=\dd P_tu+\dd\Delta^{-1/2}(I-P_t)\sqrt{\Delta}\,u .
\end{align*}
For the last equality, we use $P_t\sqrt{\Delta}\,u=\sqrt{\Delta}\,P_tu$
and apply \eqref{eq:polar} to $P_tu$.

On $S_t$, \eqref{eq:dompoisson} gives
$|\dd P_tu|^2\le C_0P_t(|\dd u|^2)<|\dd u|^2/4$.
The triangle inequality therefore gives
\[
 |\dd u|\le|\dd P_tu|+|\dd\Delta^{-1/2}(I-P_t)\sqrt{\Delta}\,u|
 \le\tfrac12|\dd u|+|\dd\Delta^{-1/2}(I-P_t)\sqrt{\Delta}\,u|.
\]
Subtracting $\tfrac12|\dd u|$ and squaring gives
$|\dd u|^2\le4|\dd\Delta^{-1/2}(I-P_t)\sqrt{\Delta}\,u|^2$ on $S_t$.
Since $\one_{S_t}=0$ on $M\setminus S_t$, we obtain
\begin{equation}\label{eq:gradient-on-stopping-set}
 0\le|\dd u|^2\one_{S_t}
 \le4|\dd\Delta^{-1/2}(I-P_t)\sqrt{\Delta}\,u|^2
 \qquad\text{a.e. on }M.
\end{equation}

To apply $P_t$ to \eqref{eq:gradient-on-stopping-set}, we check that
the right-hand side belongs to $L^1(M)$. Since $\Pi P_t=\Pi$, we have
$\Pi(I-P_t)\sqrt{\Delta}\,u=0$. Thus \eqref{eq:polar} gives
\[
 \int_M|\dd\Delta^{-1/2}(I-P_t)\sqrt{\Delta}\,u|^2\,\dif\mu
 =\|(I-P_t)\sqrt{\Delta}\,u\|_2^2
 \le\|\sqrt{\Delta}\,u\|_2^2<\infty .
\]
The inequality follows from the spectral theorem, since the multiplier
of $I-P_t$ at spectral value $\lambda$ of $\sqrt{\Delta}$ is
$1-e^{-t\lambda}\in[0,1]$. By \eqref{eq:gradient-on-stopping-set}
and $L^1$ contractivity of $P_t$, both Poisson averages are finite
almost everywhere. Integrating \eqref{eq:gradient-on-stopping-set}
against the nonnegative kernel $p_t^{\rm P}(x,y)$ in $y$ gives
\begin{align*}
 P_t\bigl(|\dd u|^2\one_{S_t}\bigr)(x)
 &=\int_M p_t^{\rm P}(x,y)|\dd u(y)|^2\one_{S_t}(y)\,\dif\mu(y)\\
 &\le4\int_M p_t^{\rm P}(x,y)
       \big|[\dd\Delta^{-1/2}(I-P_t)\sqrt{\Delta}\,u](y)\big|^2\,\dif\mu(y)\\
 &=4P_t\bigl(|\dd\Delta^{-1/2}(I-P_t)\sqrt{\Delta}\,u|^2\bigr)(x).
\end{align*}

Taking $t=2^i$ in the preceding pointwise inequality gives
\[
 P_{2^i}\bigl(|\dd u|^2\one_{S_{2^i}}\bigr)
 \le 4P_{2^i}\!\left(
   |\dd\Delta^{-1/2}(I-P_{2^i})\sqrt{\Delta}\,u|^2\right)
 \qquad\text{a.e. on }M.
\]
Since $\mathbb Z$ is countable, these inequalities hold for every $i$
outside a single null set. Taking the supremum in $i$ and then the
$L^{p/2}$ norm gives the first inequality below. The second follows
from \eqref{eq:maxsq}, since $\sqrt{\Delta}\,u\in L^2(M)\cap L^p(M)$.
\begin{align}
 Term_1
 &\le C_{p,C_0}\left\|\sup_i P_{2^i}\!\left(
   |\dd\Delta^{-1/2}(I-P_{2^i})\sqrt{\Delta}\,u|^2\right)\right\|_{p/2}
 \le C_{p,C_0}\|\sqrt{\Delta}\,u\|_p^2 .
 \label{eq:gradient-stop}
\end{align}

It remains to consider $Term_2$. If $\mu(M)<\infty$,
the definition of $\Pi$, the identity $\|\dd u\|_2=\|\sqrt{\Delta}\,u\|_2$,
and H\"older's inequality give
\[
Term_2
 =\mu(M)^{2/p-1}\|\dd u\|_2^2
 =\mu(M)^{2/p-1}\|\sqrt{\Delta}\,u\|_2^2
 \le\|\sqrt{\Delta}\,u\|_p^2 .
\]
If $\mu(M)=\infty$, then $\Pi(|\dd u|^2)=0$ and hence $Term_2=0$.
Thus, by \eqref{dd u},
\eqref{eq:gradient-stop} and the above estimates for $Term_2$, we have
\[
 \|\dd u\|_p^2
 \le C_{p,C_0}\|\sqrt{\Delta}\,u\|_p^2 .
\]
Taking square roots proves \eqref{eq:form-to-R}.

\subsection{The transform on test functions}
Let $f\in C_c^\infty(M)$. To apply \eqref{eq:form-to-R}, take
$u_b=\int_0^bP_s f\,\dif s$, $b>0$. We have
\begin{equation}\label{eq:primitive}
 \sqrt{\Delta}\,u_b=(I-P_b)f,\qquad
 \dd u_b=\dd\Delta^{-1/2}(I-P_b)f=\int_0^b\dd P_s f\,\dif s .
\end{equation}
Indeed, $f\in\operatorname{Dom}(\sqrt{\Delta})$, so $s\mapsto P_sf$
is continuous in the form norm on $[0,b]$. Thus $u_b\in W^{1,2}(M)$
and its gradient is the integral of $\dd P_s f$. At spectral value
$\lambda\ge0$ of $\sqrt{\Delta}$, the equality
$\lambda\int_0^b e^{-s\lambda}\,\dif s=1-e^{-b\lambda}$ gives the first
identity in \eqref{eq:primitive}. Applying \eqref{eq:polar} to $u_b$
gives $\dd u_b=\dd\Delta^{-1/2}(I-P_b)f$.

We check the two $L^p$ hypotheses in \eqref{eq:form-to-R}.
First, $\|\sqrt{\Delta}\,u_b\|_p=\|(I-P_b)f\|_p\le2\|f\|_p$ by
contractivity. Second, \eqref{eq:dompoisson} gives
$$\|\dd P_s f\|_p^2\le C_0\|P_s(|\dd f|^2)\|_{p/2}\le C_0\|\dd f\|_p^2.$$
Hence $$\|\dd u_b\|_p\le\int_0^b\|\dd P_s f\|_p\,\dif s
\le b\sqrt{C_0}\|\dd f\|_p<\infty.$$ We can now apply \eqref{eq:form-to-R}
to $u_b$. Using \eqref{eq:primitive} on both sides gives
\begin{align}\label{eq:uniform}
 \|\dd\Delta^{-1/2}(I-P_b)f\|_p
 &=\|\dd u_b\|_p\le C_{p,C_0}\|\sqrt{\Delta}\,u_b\|_p
 =C_{p,C_0}\|(I-P_b)f\|_p\le2C_{p,C_0}\|f\|_p .
\end{align}
As $b\to\infty$, $P_bf\to\Pi f$ in $L^2$ and $(\dd\Delta^{-1/2})\Pi=0$.
Hence $\dd\Delta^{-1/2}(I-P_b)f\to \dd\Delta^{-1/2}f$ in $L^2(TM)$. Fatou's lemma along an
almost everywhere convergent subsequence in \eqref{eq:uniform} proves
$$\|\dd\Delta^{-1/2}f\|_p\le C_{p,C_0}\|f\|_p$$ for test functions.
Density gives the bounded extension to $L^p$. On $L^2\cap L^p$ this
extension agrees with the $L^2$ realization by approximation in both norms.

For $2<p<4$, apply Riesz--Thorin interpolation to
$f\mapsto\langle \dd\Delta^{-1/2}f,e\rangle$ with any measurable complex tangent vector field
$|e|\le1$ and the inner product linear in its first variable. By
\eqref{eq:polar} and the result at $p=4$,
\[
 \|\langle \dd\Delta^{-1/2}f,e\rangle\|_p\le C_{4,C_0}^{\theta}\|f\|_p,
 \qquad \frac1p=\frac{1-\theta}{2}+\frac{\theta}{4}.
\]
For each test function $f$, take $e=(\dd\Delta^{-1/2}f)/|\dd\Delta^{-1/2}f|$ on $\{|\dd\Delta^{-1/2}f|>0\}$ and zero
elsewhere. The same estimate follows for $\|\dd\Delta^{-1/2}f\|_p$, and density gives
the assertion on $L^p$.

For $u\in C_c^\infty(M)$, the square-root formula
\[
 \sqrt{\Delta} u=\frac1{2\sqrt\pi}\int_0^\infty (u-H_su)\,\frac{\dif s}{s^{3/2}}
\]
follows from $\int_0^\infty (1-e^{-v})v^{-3/2}\,\dif v=2\sqrt\pi$
and the spectral theorem. The integral converges in $L^2$ and in $L^p$,
since $u-H_su=\int_0^s H_v\Delta u\,\dif v$ and hence
$\|u-H_su\|_r\le\min\{s\|\Delta u\|_r,2\|u\|_r\}$ for $r=2,p$.
This gives $\sqrt{\Delta} u\in L^2\cap L^p$ and $\|\dd u\|_p=\|(\dd\Delta^{-1/2})\sqrt{\Delta} u\|_p\le
C_{p,C_0}\|\sqrt{\Delta} u\|_p$.

Finally, on a disconnected manifold $\Delta$ and $\dd\Delta^{-1/2}$ act componentwise.
Applying the connected result on each component and summing the $p$th
powers gives the asserted bound with the same constant.
\hfill\qedsymbol

\bigskip
\bigskip
{\bf Acknowledgements:}
P. Chen was supported by National Key R\&D Program of China 2022YFA1005700 and by NNSF of China No. 12671127.
Duong and Li are supported by the Australian Research Council (DP260101083). L. Yan was supported by
National Key R\&D Program of China 2022YFA1005700 and by NNSF of China No. 12571111.

During preparation of the manuscript, the authors used GPT-6 for assistance with language editing and verification of our calculations. The authors take responsibility for
the mathematical arguments and the final content of the paper.

\end{document}